\documentclass{amsart}
\usepackage{amscd,amsxtra}

\newtheorem{Thm}{Theorem}
\newtheorem{Cor}{Corollary}
\newtheorem{Lem}{Lemma}
\newtheorem{Prop}{Proposition}

\theoremstyle{remark}
\newtheorem{Rem}{Remark}
\newtheorem{Def}{Definition}
\newtheorem{Ex}{Example}

\newcommand{\imbed}{\hookrightarrow}
\DeclareMathSymbol{\surj}{\mathrel}{AMSa}{"10}
\newcommand{\iso}{{\widetilde \longrightarrow}}
\newcommand{\To}{\longrightarrow}
\newcommand{\from}{\leftarrow}

\newcommand{\ot}{\otimes}
\newcommand{\oc}{\mathbin{\text{\smaller$\square$\,}}}
\DeclareMathSymbol{\birarrow}{\mathrel}{AMSa}{"13}  %rightrightarrows

\newcommand{\oplusl}{\bigoplus\limits}

\newcommand{\Hom}{{\rm Hom}}

\def\square{\hbox{\vrule\vbox{\hrule\phantom{o}\hrule}\vrule}}

\newcommand{\N}{{\mathcal N}}
\newcommand{\Nt}{{\tilde{\mathcal N}}}
\renewcommand{\O}{{\mathcal O}}
\newcommand{\F}{{\mathcal F}}
\newcommand{\G}{{\mathcal G}}

\newcommand{\Ftil}{\tilde F}

\renewcommand{\k}{{\bf k}}

\renewcommand{\b}{{\mathfrak b}}
\renewcommand{\t}{{\mathfrak t}}
\newcommand{\n}{{\mathfrak n}}
\newcommand{\g}{{\mathfrak g}}
\newcommand{\bu}{\bullet}

\newcommand{\Gm}{{\mathbb G}_m}

\newcommand{\epf}{\square}

\newcommand{\Jse}{J_\searrow}

\newcommand{\Jne}{J_\nearrow}
\newcommand{\Pnw}{P_\nwarrow}
\newcommand{\Psw}{P_\swarrow}
\newcommand{\Qsw}{Q_\swarrow}

\newcommand{\I}{{\mathbb I}}

\newcommand{\Dgood}{{D_{\infty/2}}}

\newcommand{\A}{{\mathcal A}}
\newcommand{\B}{{\mathcal B}}
\newcommand{\Zet}{{\mathbb Z}}
\newcommand{\Ce}{{\mathbb C}}
\newcommand{\cI}{{\mathcal I}}

\newcommand{\NN}{N}
\newcommand{\BB}{B}
\newcommand{\KK}{K}
\newcommand{\Nch}{N\spcheck}
\newcommand{\Ndi}{N^\#}
\newcommand{\Adi}{A^\#}
\newcommand{\phidi}{\phi^\#}
\newcommand{\Sdi}{S^\#}
\newcommand{\Bdi}{B^\#}

\newcommand{\calB}{{\mathcal B}}

\begin{document}
\title[Semi-infinite cohomology]{On semi-infinite cohomology of
 finite dimensional graded algebras}
\author{Roman Bezrukavnikov, Leonid Positselski}

\begin{abstract}
We describe a general setting for the definition of semi-infinite
cohomology of finite dimensional algebras, and provide an interpretation
of such cohomology in terms of derived categories.
%We show that semi-infinite cohomology of a finite dimensional graded
%algebra (satisfying some additional requirements) is a special case of
%a general categorical construction applied to the derived category
%of modules.
We apply this interpretation to compute semi-infinite cohomology
of some modules over the small group at a root of unity, generalizing
an earlier result of S.~Arkhipov (conjectured by B.~Feigin).
\end{abstract}
\maketitle

\begin{section}{Introduction}
Semi-infinite cohomology of associative algebras was studied
by S.~Arkhipov in \cite{Ar1}, \cite{Ar2}, \cite{Ar3}; see also \cite{S}
(these works are partly based on an earlier paper by A.~Voronov \cite{V}
where the corresponding constructions were introduced in the context
of Lie algebras).

Recall that the definition of semi-infinite cohomology
(see e.g.~\cite{Ar1}, Definition 3.3.6) works in the following set-up.
We are given an associative graded algebra $A$, two subalgebras
$\NN$, $\BB\subset A$ such that $A=\NN\otimes \BB$ as a vector space,
satisfying some additional assumptions. In this situation the space of
semi-infinite Ext's, $Ext^{\infty/2 +\bu}(X,Y)$ is defined for $X,Y$
in the appropriate derived categories.
The definition makes use of explicit complexes
(a version of the bar resolution).
The aim of this note is to show that, at least under certain simplifying
assumptions, $Ext^{\infty/2+\bu}(X,Y)$ is a particular case of
a general categorical construction.

To describe the situation in more detail, recall that starting from
an algebra $A=\NN\otimes \BB$ as above, one can define another algebra
$\Adi$, which also contains subalgebras identified with $\NN$, $\BB$,
so that $\Adi=\BB\otimes \NN$.
The semi-infinite Ext's,  $Ext^{\infty/2+\bu}(X,Y)$ are then defined
for $X\in D(\Adi-mod)$, $Y\in D(A-mod)$, where $D(\Adi-mod)$, $D(A-mod)$
are derived categories of modules with certain restrictions on
the grading.

Our categorical interpretation relies on the following construction.
Given small categories $\A$, $\A'$, $\B$ with functors  $\Phi:\B\to \A$,
$\Phi':\B\to \A'$ one can define for $X\in \A$, $Y\in \A'$ the set
of "morphisms from $X$ to $Y$ through $\B$";
we denote this set by $Hom_{\A_\B \A'}(X,Y)$.
We then show that if $\A=D^b(\Adi-mod)$, $\A'=D^b(A-mod)$, and $\B$ is
the full triangulated subcategory in $\A$ generated by $\NN$-injective
$\Adi$-modules, then, $\B$ is identified with a full subcategory
in $\A'$ generated by $\NN$-projective $A$-modules, and, under certain
assumptions, one has
\begin{equation}\label{eq}
Ext^{\infty/2 +i}(X,Y)=Hom_{\A_\B\A'}(X,Y[i]).
\end{equation}
%\footnote{Comparing our set-up with that of, say,
%\cite{Ar1}, \cite{S} the reader
% should keep in mind that the restrictions on $A,\NN$
% imposed below imply that the algebra $A^\#$
%which plays an essential role in {\it loc. cit.} is canonically
%Morita equivalent to $A$.  Thus we can talk about
%$Ext^{\infty/2 +i}(X,Y)$ for $X,Y\in D^b(A-mod)$, while in those papers
%one takes $X\in D^b(A^\#-mod)$, $Y\in D^b(A-mod)$.}

Notice that  description \eqref{eq} of $Ext^{\infty/2 +i}(X,Y)$ is
"internal" in the derived category, i.e. refers only to the derived
categories and their full subcategories rather than
 to a particular category of complexes.

\medskip

An example of the situation considered in this paper is provided by
a small quantum group at a root of unity \cite{qq},
or by the restricted enveloping algebra of a simple
Lie algebra in positive characteristic.
Computation of semi-infinite cohomology in the former case is due
to S.~Arkhipov \cite{Ar1} (the answer suggested as a conjecture by
B.~Feigin). An attempt to find a natural interpretation of this answer
was the starting point for the present work.
In section \ref{quant} we sketch a generalization of Arkhipov's Theorem
based on our description of semi-infinite cohomology and the results of
\cite{ABG}, \cite{Anichka}.
Similarly, the main result of \cite{Humph} yields a description
of semi-infinite cohomology of tilting modules over the ``big'' quantum group
restricted to the small quantum group as cohomology with support of 
coherent IC sheaves on the nilpotent cone \cite{izvrat}.

\medskip 

It should be noted that some definitions of semi-infinite
cohomology found in the literature apply in a more general (or different)
situation than the one  considered in the present paper.
An important example is provided by  affine Lie algebras; in fact,
semi-infinite cohomology has first been defined in this context,
related to the physical notion of BRST reduction.  
We hope that our approach can be extended to such  more general setting.
Some of the ingredients needed for the generalization are provided
by \cite{Posic}.

\medskip

The paper is organized as follows. 
Section \ref{sec1} is devoted to basic general facts about ``Hom through
a category''. Section 3 contains the definition of the algebra $A^\#$
and its propeties. In section 4 we recall the definition of semi-infinite
cohomology in the present context. In section 5 we prove the main result
linking that definition to the general categorical construction of
section \ref{sec1}. In section 6 we discuss the example of a small quantum
group.

\medskip

{\bf Acknowledgements.}
We are grateful to S.~Arkhipov for helpful discussions.  
This work owes its existence to W.~Soergel --- when refereeing 
the (presently unpublished) preprint \cite{preprint} submitted to
the Journal of Algebra he suggested to extend the results to
a greater generality; this is accomplished in the present paper.
We thank Wolfgang for the stimulating suggestion.  
R.B. was partially supported by
an NSF grant. He worked on this paper while visiting Princeton IAS,
the stay was funded through grants of Bell Companies, Oswald Veblen Fund,
James D Wolfenson Fund and The Ambrose Monell Foundation.
L.P. acknowledges the financial support from CRDF, INTAS, and
P.~Deligne's 2004 Balzan prize.
\end{section}

\begin{section}{Morphisms through a category}\label{sec1}

\subsection{Generalities}
Let $\A$, $\A'$, $\B$ be small categories, and $\Phi:\B\to \A$,
$\Phi':\B\to \A'$ be functors.
Fix $X\in Ob(\A)$, $Y\in Ob(\A')$.
We define the set of "morphisms from $X$ to $Y$ through
$\B$" % (or "with support at $\Phi$")
as  $\pi_0$
of  the category of diagrams
\begin{equation}\label{dia}
X\To \; \Phi(Z);\ \ \Phi'(Z)\To \; Y, \ \ \ \ \ \ \ \ \ Z\in \B.
\end{equation}
 This set will be denoted by $Hom_{\A_\B \A'}(X,Y)$.
 Thus elements of $Hom_{\A_\B \A'}(X,Y)$
are diagrams of the form \eqref{dia}, with two
diagrams identified if there exists a morphism between them.
%Composing the two arrows in \eqref{dia} we get a functorial map
%\begin{equation}\label{canmorph}
%Hom_{\A_\Phi}(X,Y)\To Hom_\A(X,Y),
%\end{equation}
%which is not injective in general.

If the categories and the functors are additive
(respectively, $R$-linear for a commutative ring $R$), then 
$Hom_{\A_\B \A'}(X,Y)$ is an abelian group (respectively, an $R$-module);
to add two diagrams of the form \eqref{dia} one sets $Z=Z_1\oplus Z_2$
with the obvious arrows.

%then addition of
%diagrams of the form \eqref{dia} is defined by
%$$(X\overset{f}{\to} \Phi( Z)\overset{g}{\to} Y)+(X
%\overset{f'}{\to} \Phi(Z')\overset{g'}{\to} Y)=
%(X\overset{f\times f'}{\To} \Phi(Z\oplus Z') \overset{g\oplus g'}{\To} Y);$$
%it induces an abelian group structure on $Hom_{\A_\Phi}(X,Y)$.
% Proposition 3 in \cite{Ma}, VIII.2 shows that  for $Z\in \B$
%the tautological map
%$$Hom(X,\Phi (Z)) \otimes _{\Zet} Hom(\Phi(Z),Y)\to Hom_{\A_\Phi}
%(X,Y)$$ is compatible with addition.

We have the composition map
$$
Hom_\A(X',X)\times Hom_{\A_\B \A'}(X,Y)\times  Hom_{\A'}(Y,Y') \to
Hom_{\A_\B\A'}(X',Y');
$$
in particular, in the additive setting
$Hom_{\A_\B\A'}(X,Y)$ is an $End(X)-End(Y)$ bimodule.

%Given $\Phi:\A\to \B$, $\Phi':\A'\to \B'$ and $F:\A\to \A'$, $G:\B\to \B'$
%with $F\circ \Phi \cong \Phi'\circ G$ we get for $X,Y\in \A$ a map
%\begin{equation}\label{F}
%Hom_{\A_\Phi}(X,Y)\to Hom_{\A'_{\Phi'}}(F(X),F(Y)).
%\end{equation}

\subsection{Pro/Ind representable case}
If the left adjoint functor $\Phi_L$ to $\Phi$ is defined on $X$,
then we have
$$Hom_{\A_\B\A'}(X,Y)=Hom_{\A'} (\Phi'(\Phi_L(X)),Y),$$
 because in this case the above category
 contracts to the subcategory of diagrams of the form 
$$X\overset{can}{\To}
\Phi(\Phi_L(X));\ \ \Phi'(\Phi_L(X)) \to Y,$$
where $can$ stands for the adjunction morphism.
 If the right adjoint functor $\Phi_R'$ is defined on $Y$, then
$$Hom_{\A_\B\A}(X,Y)=Hom_\A(X, \Phi( \Phi'_R(Y)))$$ for similar reasons.
%In particular, if $\Phi$ is a full embedding then
%\eqref{canmorph} is an isomorphism
%provided either $X$ or $Y$ lie in the image of $\Phi$.

More generally, we have

\begin{Prop}\label{proind}
Fix $X\in \A$ and $Y\in \A'$.
Assume that the functor $F_X:\B \to {\rm Sets}$,
$Z\mapsto Hom_\A(X,\Phi(Z))$ can be represented as a filtered inductive
limit of representable functors $Z\mapsto Hom_\B(\iota(S),Z)$, where
$S\in\cI$ and $\iota:\cI\to B$ is a functor between small categories.
Then we have 
$$
 Hom_{\A_\B \A'} (X,Y)=\varinjlim\limits_{S\in\cI}Hom_{A'}(\Phi'\iota(S),Y).
$$
Alternatively, assume that the functor $F_Y:\B^{op}\to {\rm Sets}$,
$Z\mapsto Hom(\Phi'(Z),Y)$ can be represented as a filtered inductive
limit of representable functors $Z\mapsto Hom_\B(Z,\iota(S))$, 
where $S\in\cI$.  Then
$$
 Hom_{\A_\B \A'} (X,Y)=\varinjlim\limits_{S\in\cI}Hom_{\A}(X,\Phi\iota(S)).
$$
\epf
\end{Prop}

\begin{Rem}\label{proindrem}
We will only use the Proposition in the case when the category $\cI$
is the ordered set of positive (or negative) integers. 
\end{Rem}

\begin{Rem}
 The assumptions of the Proposition can be rephrased by saying,
in the first case, that the functor $F_X$ is represented by
the pro-object $\varprojlim\iota$, and in the second case, that
the functor $F_Y$ is represented by the ind-object $\varinjlim\iota$.
\end{Rem}

\begin{Rem}
 The results of the Proposition can be further generalized as follows.
 Fix $X\in\A$ and $Y\in \A'$; let $\iota:\B'\to B$ be a functor
between small categories.
 Assume that either for any morphism $X\to \Phi(Z)$ the category
of pairs of morphisms $X\to\Phi\iota(S)$, $\iota(S)\to Z$ making
the triangle $X\to\Phi\iota(S)\to\Phi(Z)$ commutative is non-empty
and connected, or for any morphism $\Phi'(Z)\to Y$ the category of
pairs of morphisms $Z\to\iota(S)$, $\Phi'\iota(S)\to Y$ making
the triangle $\Phi'(Z)\to\Phi'\iota(S)\to Y$ commutative is
non-empty and connected.
 Then the natural map $Hom_{\A_{\B'} \A'} (X,Y)\to Hom_{\A_\B \A'} (X,Y)$
is an isomorphism.
\end{Rem}

\begin{Ex}\label{ex}
 Let $M$ be a Noetherian scheme, and $\A=\A'=D^b(Coh_M)$ be the
bounded derived category of coherent sheaves on $M$; let
$\Phi=\Phi':\B\imbed \A$ be the full embedding of the subcategory of
complexes whose cohomology 
is supported on a closed subset $i:N\imbed M$. Then the right adjoint functor
$i_* \circ i^!$ is well-defined as a functor to a
"larger" derived category of quasi-coherent sheaves, 
%(i.e. ind-coherent sheaves), 
while the left adjoint functor $ i_* \circ i^*$ is a
well-defined functor to the Grothendieck-Serre dual category, the
derived category of pro-coherent sheaves (introduced in Deligne's
appendix to \cite{H}). 

Let $C^\bu$ be a complex of coherent sheaves representing the object
$X\in D^b(Coh_M)$. Let $X_n$ be the 
object in the derived category represented
by the complex $C^i_n=C^i\otimes  \O_M/{\mathcal J}_N^n$
(the nonderived tensor product)
where ${\mathcal J}_N$ is the ideal sheaf of $N$. For $\F\in \B$
we have $\varinjlim Hom (X_n,\F)\iso Hom(X,\F)$. Thus
applying Proposition \ref{proind}
to $\iota:\Zet_+\to \B$ given by $\iota:n\mapsto X_n$,
we get:
 $$Hom_{\A_\B \A}(X,Y)=\varinjlim Hom(X_n,Y) =
 Hom (i_*(i^*(X)), Y)= Hom (X,i_*(i^!(Y))).$$
In particular, if $X=\O_M$ is the structure sheaf, we get
\begin{equation}\label{cohsup}
Hom_{\A_\B}(\O_M,Y[i])=H_N^i(Y),
\end{equation}
where  $H_N^\bu(Y)$ stands for cohomology with support on $N$ (see
e.g. \cite{H}).
\end{Ex}

\subsection{Triangulated full embeddings}
In all examples below $\A$, $\A'$, $\B$ will be  triangulated,
and  $\Phi$, $\Phi'$ will be full embeddings of a thick subcategory.
Assume that this is the case, and moreover $\A=\A'$, $\Phi=\Phi'$. 

\begin{Prop}\label{les}
We have a long exact sequence
$$Hom_{\A_\B \A}(X,Y) \to Hom_\A(X,Y)
 \to Hom_{\A/\B}(X,Y) \to Hom_{\A_\B \A}(X,Y[1]).$$
\end{Prop}

\proof
 The connecting homomorphism $Hom_{\A/\B}(X,Y) \to Hom_{\A_\B \A}(X,Y[1])$
is constructed as follows.
 Let $X\from X'\to Y$ be a fraction of morphisms in $\A$ representing
a morphism $X\to Y$ in $\A/\B$; the cone $K$ of the morphism $X'\to X$
belongs to $\B$.
 Assign to this fraction the diagram $X\to K$; \ $K\to Y[1]$, where
the morphism $K\to Y[1]$ is defined as the composition $K\to X'[1]\to
Y[1]$.

 All the required verifications are straightforward; the hardest one
is to check that the sequence is exact at the term $Hom_{\A/\B}(X,Y)$.
 Here one shows that for any two diagrams $X\to K'$; \ $K'\to Y[1]$ and
$X\to K''$; \ $K''\to Y[1]$ connected by a morphism $K'\to K''$ making
the two triangles commute, and for any two fractions
$X\from X'\to Y$ and $X\from X''\to Y$ to which the connecting
homomorphism assigns the respective diagrams, one can construct 
a morphism $X'\to X''$ making the triangle formed by $X'$, $X''$, $X$
commutative, and the triangle formed by $X'$, $X''$, $Y$ will then
commute up to a morphism $X\to Y$.
\epf

\end{section}

\begin{section}{Algebra $\Adi$ and modules over it}
\label{assu}

All algebras below will be associative and unital algebras over a field $\k$.

\subsection{The set-up}   \label{assumptions}

{\it We make the following assumptions.}
A $\Zet$-graded finite dimensional algebra $A$ and graded subalgebras
$\KK=A^0$, $\BB=A^{\leq 0}$, $\NN=A^{\geq 0}\subset A$ are fixed and satisfy
the following conditions:

(1) $\BB=A^{\leq 0}$, $\NN=A^{\geq 0}$ are graded by, respectively,
$\Zet^{\leq 0}$, $\Zet^{\geq 0}$, and $\KK=\BB\cap \NN$ is the
component of degree 0 in $\NN$. %and in $A^{\leq 0}$.

(2) $\KK=A^0$ is semisimple and the map $\NN\ot_\KK \BB \to A$
provided by the multiplication map is an isomorphism.

(3) Consider the $\KK$-$\NN$-bimodule $\Nch=Hom_{\KK^{op}}(\NN,\KK)$.
We require that the tensor product $S=\Nch\ot_\NN A$ is an injective 
right $\NN$-module.

\subsection{$\NN$-modules, $\Nch$-comodules, and $\Ndi$-modules}

 By a "module" we will mean a finite dimensional graded left module,
unless stated otherwise (though all the results of this section
are also applicable to ungraded or infinite dimensional modules).

 Since $\NN$ is a finitely generated projective right $\KK$-module, 
the $\KK$-bimodule $\Nch$ has a natural structure of a coring, i.e., there
is a comultiplication map $\Nch\to\Nch\ot_\KK\Nch$ and a counit map
$\Nch\to\KK$ satisfying the usual coassociativity and counity conditions.
 Consequently, there is a natural algebra structure on
$\Ndi=Hom_{\KK^{op}}(\Nch,\KK)$ and an injective morphism of algebras
$\KK\to\Ndi$.
 The category of right $\NN$-modules is isomorphic to the category of
right $\Nch$-comodules and the category of left $\Ndi$-modules is
isomorphic to the category of left $\Nch$-comodules.
 In particular, $\Nch$ is an $\Ndi$-$\NN$-bimodule.

 Recall that the cotensor product $P\oc_{\Nch}Q$ of a right $\Nch$-comodule
$P$ and a left $\Nch$-comodule $Q$ is defined as the kernel of the pair
of maps $P\ot_\KK Q\birarrow P\ot_\KK\Nch\ot_\KK Q$ one of which is induced
by the coaction map $P\to P\ot_\KK\Nch$ and the other by the coaction map
$Q\to \Nch\ot_\KK Q$.
 There are natural isomorphisms $P\oc_{\Nch}\Nch\cong P$ and $\Nch\oc_{\Nch} Q
\cong Q$.

\begin{Prop}\label{tensor-cotensor}
a) i) For any right $\NN$-module $P$ and any left $\NN$-module $Q$
there is a natural map of\/ $\k$-vector spaces
$P\ot_\NN Q \to P\oc_{\Nch}(\Nch\ot_\NN Q)$, which is an isomorphism,
at least, when $P$ is injective or $Q$ is projective.

ii) For any right $\NN$-module $P$ and any left $\Ndi$-module $Q$
there is a natural map of\/ $\k$-vector spaces
$P\ot_\NN(\NN\oc_{\Nch}Q) \to P\oc_{\Nch}Q$, which is an isomorphism,
at least, when $P$ is projective or $Q$ is injective.

b) The functors $P\mapsto\Nch\ot_\NN P$ and $M\mapsto \NN\oc_{\Nch}M$ are
mutually inverse equivalences between the categories of projective left
$\NN$-modules and injective left $\Ndi$-modules.

c) The functors $P\mapsto\Nch\ot_\NN P$ and $M\mapsto \NN\oc_{\Nch}M$ are
mutually inverse tensor equivalences between the tensor category of
$\NN$-bimodules that are projective left $\NN$-modules with
the operation of tensor product over~$\NN$ and the tensor category of
$\Ndi$-$\NN$-bimodules that are injective left $\Ndi$-modules with
the operation of cotensor product over~$\Nch$.
\end{Prop}

\proof
 Both assertions of (a) state existence of associativity (iso)morphisms
connecting the tensor and cotensor products.
 In particular, in (i) we have to construct a natural map
$(P\oc_{\Nch}\Nch)\ot_\NN Q \to P\oc_{\Nch}(\Nch\ot_\NN Q)$.
 More generally, let us consider an arbitrary $\Ndi$-$\NN$-bimodule $R$
and construct a natural map
$(P\oc_{\Nch}R)\ot_\NN Q \to P\oc_{\Nch}(R\ot_\NN Q)$.
 This map can be defined in two equivalent ways.
 The first approach is to take the tensor product of the exact
sequence of right $\NN$-modules 
$0\to P\oc_{\Nch}R\to P\ot_\KK R\to P\ot_\KK\Nch\ot_\KK R$
with the left $\NN$-module $Q$.
 Since the resulting sequence is a complex, there exists a unique map
$(P\oc_{\Nch}R)\ot_\NN Q \to P\oc_{\Nch}(R\ot_\NN Q)$ making a commutative
triangle with the natural maps of $(P\oc_{\Nch}R)\ot_\NN Q$ and
$P\oc_{\Nch}(R\ot_\NN Q)$ into $P\ot_\KK R\ot_\NN Q$.
 It is clear that this map is an isomorphism whenever $Q$ is a flat
$\NN$-module.
 Analogously, for any $P$, $Q$, $R$ there is a natural isomorphism
$(P\oc_{\Nch}R)\ot_\KK Q \cong P\oc_{\Nch}(R\ot_\KK Q)$, since $\KK$ is
semisimple.
 The second way is to take the cotensor product of the exact sequence of
left $\Nch$-comodules $R\ot_\KK\NN\ot_\KK Q\to R\ot_\KK Q\to R\ot_{\NN}Q\to0$
with the right $\Nch$-comodule $P$.
 Again, since the resulting sequence is a complex, there exists a unique
map $(P\oc_{\Nch}R)\ot_\NN Q \to P\oc_{\Nch}(R\ot_\NN Q)$ making
a commutative triangle with the natural maps from $P\oc_{\Nch}R\ot_\KK Q$
to $(P\oc_{\Nch}R)\ot_\NN Q$ and $P\oc_{\Nch}(R\ot_\NN Q)$.
 Clearly, this map is an isomorphism whenever $P$ is a coflat
$\Nch$-comodule (i.e., the cotensor product with $P$ preserves
exactness).
 Now any injective right $\NN$-module is a coflat right $\Nch$-comodule, 
since it is a direct summand of a direct sum of copies of $\Nch$.
 The two associativity maps that we have constructed coincide, since
the relevant square diagram commutes.
 The proof of (ii) is analogous.

 To prove (b), notice the isomorphisms $\NN\oc_{\Nch}(\Nch\ot_\NN P)\cong
\NN\ot_\NN P\cong P$ and $\Nch\ot_\NN (\NN\oc_{\Nch} M)\cong \Nch\oc_{\Nch}M
\cong M$ for a projective left $\NN$-module $P$ and an injective left
$\Ndi$-module $M$.
 Since a projective left $\NN$-module is a direct summand of an
$\NN$-module of the form $\NN\ot_\KK V$ and an injective left $\Ndi$-module
is a direct summand of an $\Ndi$-module of the form $\Hom_\KK(\Ndi,V)\cong
\Nch\ot_\KK V$ for a $\KK$-module $V$, the functors in question transform
projective $\NN$-modules to injective $\Ndi$-modules and vice versa.

 To deduce (c), notice the isomorphism $(\Nch\ot_\NN P)\oc_{\Nch}
(\Nch\ot_\NN Q)\cong\Nch\ot_\NN P\ot_\NN Q$ for a $\NN$-bimodule $P$
and a projective left $\NN$-module $Q$.
 It is straightforward to check that this isomorphism preserves
the associativity constraints.
\epf

\subsection{Definition of $\Adi$}

 It follows from the condition (2) that $A$ is a projective left
$\NN$-module.
 By Proposition~\ref{tensor-cotensor}(c), the tensor product 
$S=\Nch\ot_\NN A$ is a ring object in the tensor category of
$\Nch$-bicomodules with respect to the cotensor product over $\Nch$.
 By the condition (3) and the right analogue of
Proposition~\ref{tensor-cotensor}(c), the cotensor product
$\Adi=S\oc_{\Nch}\Ndi$ is a ring object in the tensor category of
$\Ndi$-bimodules with respect to the tensor product over $\Ndi$.
 The embedding $\NN\to A$ induces injective maps $\Nch\to S$ and
$\Ndi\to\Adi$; these are unit morphisms of the ring objects in
the corresponding tensor categories.
 So $\Adi$ has a natural associative algebra structure and $\Ndi$
is identified with a subalgebra in $\Adi$.
 Notice that $\Adi$ is a projective right $\Ndi$-module by
the definition.

\begin{Prop} \label{s-sdi}
 There is a natural isomorphism between the $\Ndi$-$A$-bimodule
$S=\Nch\ot_\NN A$ and the $\Adi$-$\NN$-bimodule $\Sdi=\Adi\ot_{\Ndi}\Nch$,
making $S$ an $\Adi$-$A$-bimodule.
 Moreover, there are isomorphisms:
$$
 \Adi\cong End_{A^{op}}(S), \quad A^{op}\cong End_{\Adi}(\Sdi).
$$
\end{Prop}

\proof
 By the definition, we have $\Sdi=(S\oc_{\Nch}\Ndi)\ot_{\Ndi}\Nch\cong
S\oc_{\Nch}\Nch\cong S$, since $S$ is an injective right $\NN$-module.
 Let us show that the right $A$-module and the left $\Adi$-module
structures on $S\cong\Sdi$ commute.
 The isomorphism $S\ot_\NN A\cong S\ot_\NN(\NN\oc_{\Nch}S)\cong S\oc_{\Nch}S$
transforms the right action map $S\ot_\NN A\to S$ into the map
$S\oc_{\Nch}S\to S$ defining the structure of ring object in the tensor
category of $\Nch$-bicomodules on $S$.
 Analogously, the isomorphisms $\Adi\ot_{\Ndi}\Sdi\cong(S\oc_{\Nch}\Ndi)
\ot_{\Ndi}S\cong S\oc_{\Nch} S$ and $\Sdi\cong S$ transform the left action
map $\Adi\ot_{\Ndi}\Sdi\to\Sdi$ into the same map $S\oc_{\Nch}S\to S$.
 Finally, there is an isomorphism $\Adi\ot_{\Ndi}S\ot_\NN A\cong
(S\oc_{\Nch}\Ndi)\ot_{\Ndi}S\ot_\NN(\NN\oc_{\Nch}S)\cong
S\oc_{\Nch}S\oc_{\Nch}S$, so the right and left actions commute since
$S$ is an associative ring object in the tensor category of
$\Nch$-bicomodules.
 Now we have $Hom_{A^{op}}(\Nch\ot_N A, \Nch\ot_N A)\cong
Hom_{N^{op}}(\Nch, \Nch\ot_N A)\cong (\Nch\ot_N A)\oc_{\Nch}\Ndi = \Adi$
and $Hom_{\Adi}(\Adi\ot_{\Ndi}\Nch, \Adi\ot_{\Ndi}\Nch)\cong
Hom_{\Ndi}(\Nch, \Adi\ot_{\Ndi}\Nch)\cong N\oc_{\Nch}(\Adi\ot_{\Ndi}\Nch)
\cong A$.
\epf

\subsection{$\NN$-projective (injective) modules}
 By $A-mod$ we denote the category of (graded finite dimensional)
left $A$-modules.

Consider the full subcategories $A-mod_{\NN-proj}\subset A-mod$,
$\Adi-mod_{\Ndi-inj} \subset \Adi-mod$ consisting of modules whose
restriction to $\NN$ is projective (respectively, restriction
to $\Ndi$ is injective).

We abbreviate $D(A)=D^b(A-mod)$, $D(\Adi)=D^b(\Adi-mod)$, and let
$D_{\infty/2}(A)\subset D(A)$, $D_{\infty/2}(\Adi)\subset D(\Adi)$ be
the full triangulated subcategories generated by $A-mod_{\NN-proj}$,
$\Adi-mod_{\Ndi-inj}$ respectively.

\begin{Prop}
We have canonical equivalences:  $A-mod_{\NN-proj}\cong \Adi-mod_{\Ndi-inj}$,
$D_{\infty/2}(A) \cong D_{\infty/2}(\Adi)$.
\end{Prop}

\proof
 Let us show that the adjoint functors $P\mapsto S\ot_A P$ and
$M\mapsto Hom_{\Adi}(S,M)$ between the categories $A-mod$ and $\Adi-mod$
induce an equivalence between their full subcategories $A-mod_{\NN-proj}$
and $\Adi-mod_{\Ndi-inj}$.
 It suffices to check that the adjunction morphisms 
$P\to Hom_{\Adi}(S,\;S\ot_A P)$ and $S\ot_A Hom_{\Adi}(S,M)\to M$ are
isomorphisms when an $A$-module $P$ is projective over $\NN$ and
an $\Adi$-module $M$ is injective over $\Ndi$.
 There are natural isomorphisms $S\ot_A P\cong \Nch\ot_N P$ and
$Hom_{\Adi}(S,M)\cong Hom_{\Ndi}(\Nch,M)\cong \NN\oc_{\Nch} M$, so it
remains to apply Proposition~\ref{tensor-cotensor}(b).
 To obtain the equivalence of categories
$D_{\infty/2}(A) \cong D_{\infty/2}(\Adi)$, it suffices to check that
$D_{\infty/2}(A)$ is equivalent to the bounded derived category of the exact
category $A-mod_{\NN-proj}$ and $D_{\infty/2}(\Adi)$ is equivalent to
the bounded derived category of the exact category $\Adi-mod_{\Ndi-inj}$.
 Let us prove the former; the proof of the latter is analogous.
 It suffices to check that for any bounded complex of $\NN$-projective
$A$-modules $P$ and any bounded complex of $A$-modules $X$ together
with a quasi-isomorfism $X\to P$ there exists a bounded complex of 
$\NN$-projective $A$-modules $Q$ together with a quasi-isomorphism
$Q\to X$.
 Let $Q'$ be a bounded above complex of projective $A$-modules mapping
quasi-isomorphically into $X$; then the canonical truncation
$Q'_{\ge -n}$ for large enough $n$ provides the desired complex $Q$.
\epf

\subsection{The case of an invertible entwining map}  \label{entwining}
Consider the multiplication map $\phi:\BB\ot_{\KK}\NN  \to A\cong 
\NN\ot_{\KK}\BB$.  It yields a map $\psi:\Nch\ot_{\KK}\BB \to
Hom_{\KK^{op}}(\NN,\BB)\cong\BB\ot_\KK\Nch$.
Assume that the map $\psi$ is an isomorphism and consider
the inverse map $\psi^{-1}: \BB\ot_{\KK}\Nch \to \Nch\ot_{\KK}\BB$.
By the analogous "lowering of indices" we obtain from it a map
$\Ndi\ot_{\KK}\BB \to Hom_{\KK^{op}}(\Nch,\BB)=\BB\ot_{\KK}\Ndi$ that will be
denoted by $\phidi$.

Then the algebra $\Adi$ can be also defined as the unique associative
algebra with fixed embeddings of $\Ndi$ and $\BB$ into $\Adi$ such that

i) the embeddings $\Ndi\to\Adi$ and $\BB\to\Adi$ form a commutative
square with the embeddings $\KK\to\Ndi$ and $\KK\to\BB$;

ii) the multiplication map induces an isomorphism
$\BB\ot_{\KK}\Ndi \to \Adi$;

iii) the map induced by the multiplication map $\Ndi\ot_{\KK}\BB \to \Adi
\cong\BB\ot_{\KK}\Ndi$ coincides with $\phidi$.

Indeed, the existence of an algebra $A$ with subalgebras $\NN$ and $\BB$
in terms of which the map $\phi$ is defined can be easily seen to
be equivalent to the map $\psi$ satisfying the equations of a right
entwining structure for the coring $\Nch$ and the algebra $\BB$
(see \cite{cobook} or \cite{Posic} for the definition).
When $\psi$ is invertible, it is a right entwining structure if and
only if $\psi^{-1}$ is a left entwining structure, and the latter is
equivalent to the existence of an algebra $\Adi$ satisfying (i-iii).

To show that the two definitions of $\Adi$ are equivalent, it
suffices to check that the ring object $S$ in the tensor category of
$\Nch$-bicomodules can be constructed in terms of the entwining
structure $\psi$ in the way explained in~\cite{copaper} or~\cite{Posic}.

\subsection{The case of a self-injective $\NN$}
Assume that $\NN$ is self-injective. 
In this case $\Adi$ is canonically Morita equivalent to $A$;
the equivalence is defined by the $\Adi$-$A$-bimodule $S$, so it sends
$A-mod_{N-proj}=A-mod_{N-inj}$ to $\Adi-mod_{\Ndi-proj}=\Adi-mod_{\Ndi-inj}$.

Indeed, $\Nch$ is obviously an injective generator of the category of
right $\NN$-modules.
Since every injective $\NN$-module is projective, $\Nch$ is a projective
right $\NN$-module.
Since $\NN$ is an injective right $\NN$-module, it is a direct summand
of a finite direct sum of copies of $\Nch$.
So $\Nch$ is a projective generator of the category of right
$\NN$-modules; hence $S=\Nch\ot_\NN A$ is a projective generator of
the category of right $A$-modules.
Now it remains to use Proposition~\ref{s-sdi}.
Analogously, $\Ndi$ is Morita equivalent to $\NN$; hence $\Ndi$ is
also self-injective.

If $\NN$ is Frobenius, $\Ndi$ is isomorphic to $\NN$ and $\Adi$ is
isomorphic to $A$.
Indeed, $\KK$ is also Frobenius; choose a Frobenius linear function
$\KK\to\k$; then the right $\KK$-module $\Hom_\k(\KK,\k)$ is isomorphic
to $\KK$.
Hence a Frobenius linear function $\NN\to\k$ lifts to a right
$\KK$-module map $\NN\to\KK$.
Now the composition $\NN\ot_\KK\NN\to \NN\to \KK$ of the multiplication
map $\NN\ot_\KK\NN\to \NN$ and the right $\KK$-module map $\NN\to\KK$
defines an isomorphism of right $\NN$-modules $\NN\to Hom_{\KK^{op}}(\NN,\KK)
=\Nch$.
By Proposition~\ref{s-sdi}, this leads to the isomorphism $\Adi\cong A$
and analogously to the isomorphism $\Ndi\cong\NN$; these isomorphisms
are compatible with the embeddings $\NN\to A$ and $\Ndi\to\Adi$, but
not with the embeddings of $\KK$ to $\NN$ and $\Ndi$, in general.

%Then the condition (3) above amounts to the requirement that
%the multiplication map induces an isomorphism $\BB\otimes \NN \to A$.

\end{section}

\begin{section}{Definitions of $Ext^{\infty/2}$ by explicit complexes}
\subsection{Concave and convex resolutions}
A complex of graded modules will be called {\it convex} if the grading
"goes down", i.e. for any $n\in \Zet$ the sum of graded components of
degree more than $n$ is finite dimensional; it will be called
{\it non-strictly convex} if the grading "does not go up", i.e.
the graded components of high enough degree vanish.
A complex of graded modules will be called {\it concave} (respectively
{\it non-strictly concave}) if the grading "goes up" (respectively
"does not go down") in the similar sense.

An $\Adi$-module $M$ will be called {\it weakly projective relative
to $\Ndi$} if for any $\Adi$-module $J$ which is injective as
an $\Ndi$-module one has $Ext_{\Adi}^i(M,J)=0$ for all $i\ne 0$.
Analogously one defines $A$-modules {\it weakly injective
relative to $\NN$.}
Notice that any $\Adi$-module induced from an $\Ndi$-module is weakly
projective relative to $\Ndi$.
The class of $\Adi$-modules weakly projective relative to $\Ndi$ is
closed under extensions and kernels of surjective morphisms.

\begin{Lem}\label{resolve_b}
i) Any $A$-module admits a left concave resolution by $A$-modules
which are projective as $\NN$-modules.
 Any $\Adi$-module admits a left non-strictly convex resolution
by $\Adi$-modules which are weakly projective relative to $\Ndi$.

ii) Any finite complex of $A$-modules is a quasiisomorphic quotient
of a bounded above concave complex of $\NN$-projective $A$-modules. 
 Any finite complex of $\Adi$-modules is a quasiisomorphic quotient
of a bounded above non-strictly convex complex of $\Adi$-modules
weakly projective relative to $\Ndi$.
\end{Lem}

\proof To deduce (ii) from (i) choose a quasiisomorphic surjection
onto a given complex $C^\bu\in Com^b(A-mod)$ from a complex of
$A$-projective modules $P^\bu\in Com^-(A-mod)$ (notice that condition
(2) of~\ref{assumptions} implies that an $A$-projective module is also
$\NN$-projective), and apply (i) to the module of cocycles
$Z^n =P^{-n}/d(P^{-n-1})$ for large $n$.

To check (i) it suffices to find for any $M\in A-mod$ a surjection
$P\surj M$, where $P$ is $\NN$-projective, and if $n$ is such that
all graded components $M_i$ for $i<n$ vanish, then $P_i=0$ for $i<n$
and $P_n\iso M_n$. 
It suffices to take $P=Ind_{\BB}^A (Res_{\BB}^A(M))$.
It is indeed $\NN$-projective, because of the equality
\begin{equation}\label{resind}
Res^A_{\NN}(Ind_{\BB}^A(M))=Ind_{\KK}^{\NN}(M)),
\end{equation}
which is a consequence of assumption (2). \epf

The second assertions of (i) and (ii) are proven in the analogous
way, except that one uses the induction from $\Ndi$ (this is even
simpler, as weak relative projectivity of the relevant modules
is just obvious).

\subsection{Definition of semi-infinite Ext's}
\begin{Def}\label{defpol} (cf. \cite{LNM}, \S 2.4)
The assumptions (1--3) of~\ref{assumptions} are enforced.
Let $X\in D(\Adi)$ and $Y\in D(A)$.
Let $\Psw^X$ be a non-strictly convex bounded above complex of
$\Adi$-modules weakly projective relative to $\Ndi$ that is
quasiisomorphic to $X$, and $\Pnw ^Y$ be a concave bounded above
complex of $\NN$-injective $A$-modules that is quasiisomorphic
to $Y$. Then we set
\begin{equation}\label{defini}
Ext^{\infty/2+i}(X,Y)=H^i(Hom^\bu_{\Adi} (\Psw^X, S\ot_A\Pnw^Y)).
\end{equation}
Independence of  the right-hand side of \eqref{defini} on the
choice of  $\Psw^X$, $\Pnw^Y$ follows from Theorem \ref{teorem}
below.
\end{Def}

\begin{Rem}\label{grad}
Notice that $Hom$ in the right-hand side of \eqref{defini} is $Hom$
in the category of graded modules. As usual, it is often convenient
to denote by $Ext^{\infty/2+i}(X,Y)$
the graded space which in present notations
is written down as $\oplusl_n Ext^{\infty/2+i}(X,Y(n))$, where
$(n)$ refers to the shift of grading by $-n$.
\end{Rem}

\begin{Rem}
Definition \ref{defpol} is compatible with \cite{Ar1}, Definition 3.3.6
in the sense explained below.
In this remark we will freely use the notation of {\it loc. cit.}

For a finite dimensional algebra $A$ the definition of the algebra
$\Adi$ given in \cite{Ar1}, 3.3.2 reduces to $\Adi=End_{A^{op}}(S)$,
where $S$ is defined by $S=Hom_k(\NN,k)\ot_{\NN}A$, so according to
Proposition~\ref{s-sdi} this agrees with our definition
(see also~\ref{entwining}).
Notice that in {\it loc. cit.} it is presumed that $K=k$, so one has
$\Ndi=\NN$.

Let $L\in Com^b(\Adi-mod)$ and $M\in Com^b(A-mod)$.
Then the restricted Bar-resolution ${\rm Bar}^\bu(\Adi,\Ndi,L)$ is
a non-strictly convex bounded above resolution of $L$ by $\Adi$-modules
weakly projective relative to $\Ndi$; and ${\rm Bar}^\bu(A,B,M)$ is 
a concave bounded above resolution of $M$ by $\NN$-projective
$A$-modules.
Thus the definition of semi-infinite cohomology
$$
 Ext^{\infty/2+i}(L,M) = Hom^\bu_{\Adi} \left(                                
 {\rm Bar}^\bu (\Adi,\Ndi,L), S\ot_A {\rm Bar}^\bu (A,B,M) \right)
$$
from {\it loc. cit.} is a particular case of our definition
whenever both are applicable.
\end{Rem}

\subsection{Alternative assumptions}
%\begin{Rem}
The conditions on the resolutions $\Psw^X$, $\Pnw^Y$ used in
\eqref{defini} are formulated in terms of the subalgebras $\NN\subset A$
and $\Ndi\subset\Adi$; the subalgebra $\BB\subset A$ is not mentioned
there (and the left-hand side of \eqref{teoreq} in Theorem \ref{teorem}
below does not depend on it either).
However, existence of a "complemental" subalgebra $\BB$ is used in
the construction of a resolution $\Pnw^Y$ with required properties.
Moreover, the next standard Lemma shows that conditions on
the resolutions $\Psw^X$, $\Pnw^Y$ can be alternatively rephrased
in terms of the subalgebra $\BB$ and any nonpositively graded
subalgebra $\Bdi\subset\Adi$ such that $\Bdi\ot_\KK\Ndi\cong\Adi$,
assuming that such a subalgebra exists (e.g., in the assumptions
of~\ref{entwining} or when $\NN$ is Frobenius and $B\ot_\KK\NN\cong A$).
%\end{Rem}

\begin{Lem}\label{fifi}
i) An $A$-module is $\NN$-projective iff it has a filtration with
subquotients of the form $Ind_\BB^A(M)$, $M\in \BB-mod$.

ii) Assume that $\Bdi\subset\Adi$ is a subalgebra graded by nonpositive
integers such that $\KK\subset\Bdi$ and the multiplication map induces
an isomorphism $\Bdi\ot_\KK\Ndi\to\Adi$.
 Then an $\Adi$-module is $\Ndi$-injective iff it has a filtration
with subquotients of the form $CoInd_{\Bdi}^{\Adi}(M)$, $M\in \Bdi-mod$.
\end{Lem}

\proof The "if" direction follows from semisimplicity of $\KK$,
and equality \eqref{resind} above.
To show the "only if" part let $M$ be a projective $\NN$-module.
Let $M^-$ be its graded component of minimal degree;
then the canonical morphism
\begin{equation}\label{mm}
Ind_\KK^\NN M^-\to M
\end{equation}
is injective. If $M$ is actually an $A$-module, then the injection
$M^-\to M$ is an embedding of $\BB$-modules, hence yields
a morphism of $A$-modules
\begin{equation}\label{mmm}
Ind_\BB^A M^-\to M.
\end{equation}
 \eqref{resind} shows that $Res^A_\NN$ sends \eqref{mmm}
into \eqref{mm}; in particular \eqref{mmm} is injective. Thus
the bottom submodule of the required filtration is constructed,
and the proof is finished by induction. 
The proof of~(ii) is analogous. \epf

\begin{Rem}  \label{injective}
 Replacing the assumption of existence of a subalgebra $\BB\subset A$
(assuming only that $A$ is a projective left $\NN$-module) with
the assumption of existence of a nonpositively graded subalgebra
$\Bdi\subset\Adi$ such that $\Bdi\ot_\KK\Ndi\cong\Adi$, one can
define $Ext^{\infty/2+i}(X,Y)$ in terms of injective resolutions rather
than projective ones.
 Namely, for $X\in D(\Adi)$ and $Y\in D(A)$, let $\Jse^X$ be
a convex bounded below complex of $\Ndi$-injective modules
quasiisomorphic to $X$, and $\Jne ^Y$ be a non-strictly concave
bounded below complex of $A$-modules weakly injective relative to
$\NN$. Then set
$$
 Ext^{\infty/2+i}(X,Y)=H^i(Hom^\bu (Hom_{\Adi}(S,\Jse ^X), \Jne^Y)).
$$
The analogue of Theorem \ref{teorem} below holds for this definition
as well, hence it follows that the two definitions are equivalent
whenever both are applicable.
\end{Rem}

%\begin{Rem} 
\subsection{Comparison with ordinary Ext and Tor}\label{ord}
In four special cases $Ext^{\infty/2+i}(X,Y)$ coincides with
a combination of traditional derived functors.
First, suppose that $Res ^A_\NN(Y)$ has finite projective dimension;
then one can use a finite complex $\Pnw^Y$ in \eqref{defini} above.
It follows immediately, that in this case we have
$$Ext^{\infty/2+i}(X,Y)\cong Hom_{D(\Adi)}(X,S \overset{L}{\otimes}_AY[i]).$$
Analogously, in the assumptions of Remark~\ref{injective} above,
whenever $Res^{\Adi}_{\Ndi}(X)$ has finite injective dimension one has
$$Ext^{\infty/2+i}(X,Y)\cong Hom_{D(A)}(RHom_{\Adi}(S,X),Y[i]).$$
On the other hand, suppose that the complex $\Psw^X$ in \eqref{defini}
can be chosen to be a finite complex of $\Adi$-modules whose terms
have filtrations with subquotients being $\Adi$-modules induced
from $\Ndi$-modules.
We claim that in this case we have
$$Ext^{\infty/2+i}(X,Y)\cong H^i(RHom_{\Adi}(X,S)\overset{L}{\otimes}_A Y).$$
This isomorphism is an immediate consequence of the next Lemma.
Analogously, in the situation of Remark~\ref{injective}, whenever
$\Jne^Y$ can be chosen to be a finite complex of $A$-modules whose terms
have filtrations with subquotients being $A$-modules coinduced from
$\NN$-modules, one has
$$Ext^{\infty/2+i}(X,Y)\cong H^i(X^* \overset{L}{\otimes}_{\Adi}
RHom_A(S^*,Y)).$$
Here we denote by $V\mapsto V^*$ the passage to the dual vector space,
$V^*=Hom_k(V,k)$, and the corresponding functor on the level of derived
categories.
%\end{Rem}

\begin{Lem} Let  $L\in \Adi-mod$, $M\in A-mod$ be such that $L$ has
a filtration with subquotients being $\Adi$-modules induced from
$\Ndi$-modules, while $M$ is $\NN$-projective. Then we have

a) i) $Ext^i_{\Adi}(L,S)=0$ and $Tor_i^A(Hom_{\Adi}(L,S),M)=0$
for $i\ne 0$.

ii) $Tor_i^A(S,M)=0$ and $Ext^i_{\Adi}(L,S\ot_A M)=0$ for $i\ne 0$.

b) The natural map
\begin{equation}\label{isomom}
Hom_{\Adi}(L,S)\ot_A M \To Hom_{\Adi}(L,S\ot_A M)
\end{equation}
is an isomorphism.
\end{Lem}

\proof The first equality in (i) holds because $S$ is an injective
$\Ndi$-module.
To check the second one, notice that if $L=Ind_{\Ndi}^AL_0$, then
$Hom_{\Adi}(L,S) \cong Hom_{\Ndi}(L_0,\Nch\ot_\NN A)\cong
Hom_{\Ndi}(L_0,\Nch)\ot_\NN A$ is a right $A$-module induced from
a right $\NN$-module.
The first equality in (ii) holds because the right $A$-module $S$
is induced from a right $\NN$-module, and the second one is verified
since $S\ot_A M$ is $\Ndi$-injective.
Let us now deduce (b) from (a). Notice that (a) implies that both
sides of \eqref{isomom} are exact in $M$ (and also in $L$),
i.e. send exact sequences $0\to M' \to M\to M'' \to 0$
with $M'$, $M''$ being $\NN$-projective into exact sequences.
Also \eqref{isomom} is evidently an isomorphism for $M=A$.
For any $\NN$-projective $M$ there exists an exact sequence
$$A^n\overset{\phi}{\To}A^m\to M\to0$$
with the image and kernel of $\phi$ being $\NN$-projective.
Thus both sides of \eqref{isomom} turn into exact sequences, which
shows that \eqref{isomom} is an isomorphism for any $\NN$-projective $M$.
\epf

\end{section}

\begin{section}{Main result}

\begin{Thm}\label{teorem}
 Let $\Dgood\subset D(\Adi)$, $\Dgood\subset D(A)$ be the full
triangulated subcategory of $D(\Adi)$ generated by $\Ndi$-injective
modules, which is equivalent to the full triangulated subcategory
of $D(A)$ generated by $\NN$-projective modules.
For
 $X\in D^b(\Adi-mod)$, $Y\in D^b(A-mod)$ we have a natural
isomorphism
\begin{equation}\label{teoreq}
Hom_{{D(\Adi)}_{\Dgood}D(A)}(X,Y[i])\cong Ext^{\infty/2+i}(X,Y).
\end{equation}
\end{Thm}

\begin{Ex}\footnote{We thank A.~Beilinson who suggested to us
this example.}
Assume that $A=N$ is a Frobenius algebra and $K=\k$. Then $\Adi\cong A$,
and according to section \ref{ord} we have $Ext^{\infty/2+i}(X,Y)=
Tor_{-i}^A(X^*,Y)$. 
In this case we can identify $\A'$ with $\A$, so that
 $\Phi'=\Phi$ is the 
embedding of the category of perfect complexes.
The long exact sequence of Proposition \ref{les} becomes a standard sequence
linking Ext, Tor and Hom in the stable category $A/\B$; in particular,
for modules over a finite group we recover the description
of Tate cohomology as Hom's in the stable category.

\end{Ex}

\begin{Rem}
Notice that the definition of the left hand side in \eqref{teoreq} 
applies also to non-graded algebras and modules.
Thus the Theorem allows one to extend the definition of semi-infinite
cohomology to nongraded algebras.
Another definition of the semi-infinite cohomology of nongraded
algebras was given in~\cite{Posic}. 

Let us point out that these two definitions are \emph{not} equivalent:
for example, when $\k$ is a finite or a countable field, the left hand
side of \eqref{teoreq} in the nongraded case is no more than countable,
while the semi-infinite cohomology defined in {\it loc. cit.} can have
the cardinality of continuum.
%as they in fact differ even for a graded algebra considered without
%the grading.
%Namely, let $A$ be an algebra satisfying (1-3) of~\ref{assumptions}
%and $X$ and $Y$ be finite complexes of graded modules over $\Adi$ and
%$A$, respectively.
%Then the semi-infinite Ext between $X$ and $Y$ defined as the left hand
%side of \eqref{teoreq} in terms of categories of nongraded modules is
%isomorphic to the direct sum of $Ext^{\infty/2+\bu}(X,Y(n))$ over
%$n\in\Zet$ (in the notation of Remark~\ref{grad}), while the nongraded
%$SemiExt_S(X,Y)$ is isomorphic to the direct product of
%the same components.
\end{Rem}

\medskip
The proof of  Theorem \ref{teorem} is based on the following

\begin{Lem}\label{resolve_a}
i) Every $\NN$-projective $A$-module admits a non-strictly convex
left resolution consisting of $A$-projective modules.

ii) A finite complex of $\NN$-projective $A$-modules is quasiisomorphic
to a non-strictly convex bounded above complex of $A$-projective modules.

\end{Lem}
\proof (ii) follows from (i) as in the proof of Lemma \ref{resolve_b}.
(Recall that, according to a well-known argument due to Hilbert,
if a bounded above complex of projectives represents an object of
the derived category which has finite projective dimension, then for
large negative $n$ the module of cocycles is projective.)

To prove (i) it is enough for any $\NN$-projective module
$M$ to find a surjection $Q\surj M$, where $Q$ is $A$-projective,
and $Q_n=0$ for $i>n$ provided $M_i=0$ for $i>n$. (Notice that
the kernel of such a surjection is  $\NN$-projective, because
$Q$ is  $\NN$-projective by condition (2).) We can take $Q$
to be $Ind_\NN^A (Res_\NN^A(M))$, and the condition on grading 
is clearly satisfied.
\epf

\begin{Prop} a) Let $\Pnw$ be a concave bounded above complex of
$A$-modules representing an object $Y\in D^-(A-mod)$.
Let $\Pnw^n$ be the $(-n)$-th stupid truncation of $\Pnw$
(thus $\Pnw^n$ is a subcomplex of $\Pnw$).

Let $Z$ be a finite complex of $\NN$-projective $A$-modules.
Then we have
\begin{equation}\label{lim}
Hom_{D^-(A-mod)}(Z,Y) \iso \varinjlim Hom_{D(A)} (Z,\Pnw^n).
\end{equation}
In fact, for $n$ large enough we have
$$Hom_{D^-(A-mod)}(Z,Y) \iso Hom_{D(A)}(Z,\Pnw^n).$$
\end{Prop}

\proof Let $\Qsw$ be a non-strictly convex bounded above complex of
$A$-projective modules quasiisomorphic to $Z$ (which exists by
Lemma \ref{resolve_a}(ii)).
Then the left-hand side of \eqref{lim} equals $Hom_{Hot}(\Qsw,\Pnw)$,
where $Hot$ stands for the homotopy category of complexes of $A$-modules.
The conditions on gradings of our complexes ensure that there are only
finitely many degrees for which the corresponding graded components both
in $\Qsw$ and $\Pnw$ are nonzero; thus any morphism between the graded
vector spaces $\Qsw$, $\Pnw$ factors through the finite dimensional sum
of the corresponding graded components. In particular,
$Hom^\bu(\Qsw, \Pnw^n) \iso Hom^\bu (\Qsw, \Pnw)$ for large $n$,
and hence $$Hom_{D(A)} (Z, \Pnw^n) = Hom_{Hot} (\Qsw, \Pnw^n) \iso
Hom_{Hot} (\Qsw, \Pnw)$$ for large $n$.
\epf

\begin{Cor}
Let $\Pnw$ be a concave bounded above complex of $\NN$-projective
$A$-modules, and $X$ be the corresponding object of $D^-(A-mod)$.
Then the functor on $\Dgood$ given by $Z\mapsto Hom_{D^-(A-mod)}(Z,Y)$
is represented by the ind-object\/  $\varinjlim \Pnw^n$. \epf
\end{Cor}

{\it Proof} of the Theorem. We keep the notation of Definition
\ref{defpol}.
It follows from the Proposition that
$$      % \begin{equation}\label{last1}
Hom_{D(\Adi)_\Dgood D(A)}(X,Y[i]) =
\varinjlim_n Hom_{D(\Adi)}(X, S\ot_A(\Pnw^Y)^n).
$$     %\end{equation}
The right-hand side of \eqref{teoreq} (defined in \eqref{defini})
equals $H^i( Hom^\bu_{\Adi}(\Psw^X, S\ot_A\Pnw^Y))$. 
The conditions on gradings of $\Psw^X$, $\Pnw^Y$ show that
for large $n$ we have
$$   %\begin{equation}\label{last}
 Hom^\bu_{\Adi}(\Psw^X, S\ot_A(\Pnw^Y)^n) \iso
Hom^\bu_{\Adi}(\Psw^X, S\ot_A\Pnw^Y) .
$$     %% \end{equation}
Since $Ext^i_{\Adi}(L,S\ot_A M)=0$ for $i>0$ if $L$ is weakly projective
relative to $\Ndi$ and $M$ is $\NN$-projective, we have
$$Hom_{D(\Adi)}(X, S\ot_A(\Pnw^Y)^n)=Hom^\bu_{\Adi}(\Psw^X, S\ot_A(\Pnw^Y)^n).$$
The Theorem is proved. \epf

\begin{Rem}
 There is a version of Theorem~\ref{teorem} applicable in the situation
when the condition that $\KK$ is the component of degree~0 of $\NN$ 
in (1) of \ref{assumptions} is replaced with the condition that $\KK$
is the component of degree~0 of $\BB$.
 One just has to change the conditions on the complexes $\Psw^X$, $\Pnw^Y$
in Definition~\ref{defpol}, requiring that $\Psw^X$ be convex and $\Pnw^Y$
be non-strictly concave, and make the related changes in the proof.
\end{Rem}

\section{Semi-infinite cohomology of the small quantum group}\label{quant}

This section concerns with the example provided by a
small quantum group. So let $\g$ be a simple Lie algebra over $\Ce$
with a fixed triangular decomposition $\g=
\n\oplus\t\oplus \n^-$. %=\n \oplus \b^-$.
Let $q\in \Ce$ be a root of unity of order $l$, and
 let $A=u_q=u_q(\g)$ be the corresponding small quantum group \cite{qq}.
We assume  that $l$ is large enough (larger than Coxeter number)  %CHECK!
and is prime to twice
the maximal multiplicity of an edge in the Dynkin
diagram of $\g$.

Let $A^{\geq 0}= u_q^+\subset u_q$ and
$A^{\leq 0}=u_q^-\subset u_q$ be respectively
 the upper and the lower triangular subalgebras.
The algebra $u_q$ carries a canonical grading by the weight lattice.
We fix an arbitrary element in the dual coweight lattice, which is a dominant
coweight, thus we obtain a $\Zet$-grading on $u_q$.
Then the above conditions (1--3) are satisfied.

For an augmented $\k$ algebra $R$ we write $H^\bu(R)$ for $Ext_R(\k,\k)$;
we abbreviate
 $H^\bu=H^\bu(u_q)$.

 The cohomology algebra $H^\bu$, and the semi-infinite cohomology
$Ext^{\infty/2+\bu}(\k,\k)$ were computed respectively in \cite{GK} and
\cite{Ar1}. Let us recall the results of these computations.
Below by ``Hom'' we will mean graded Hom  as in Remark \ref{grad}
above.

Let $\N\subset \g$ be the cone of nilpotent elements, and $\n\subset \N$ be
a maximal nilpotent subalgebra.
 Then the Theorem  of Ginzburg and Kumar
 asserts the existence of canonical isomorphisms
\begin{equation}\label{GinK}
H^\bu\cong \O(\N),
\end{equation}
\begin{equation}\label{GinKb}
H^\bu (u_q^+)=\O(\n),
\end{equation}
such that the restriction map $\O(\N)\to \O(\n)$ coincides
with the map arising from functoriality of cohomology with respect to
maps of augmented algebras.

Also,
 a Theorem of Arkhipov (conjectured by Feigin) asserts that
\begin{equation}\label{Arki}
Ext^{\infty/2+\bu}(\I,\I)\cong H_{\n^-}^d(\N,\O),
\end{equation}
where $d$ is the dimension of $\n^-$, and $H_{\n^-}$ denotes cohomology with
support in $\n^-$; one also has $H^i_{\n^-}(\N,\O)=0$ for $i\ne d$.

The aim of this section is to show how (a generalization of) this isomorphism
follows from Theorem \ref{teorem}.

\subsection{$D_{\infty/2}$ and cohomological support}

Let $D^\bu$ denote the category  defined
by $Ob(D^\bu)=Ob(D)$, $Hom_{D^\bu}(X,Y)=Hom^\bu(X,Y)=\oplusl_i Hom(X,Y[i])$.
Then $D^\bu$ is an $HH^\bu$-linear category, i.e. we have a canonical
homomorphism $HH^\bu\to End(Id_{D^\bu})$, where $HH^\bu$ denotes the
Hochschield cohomology of $u_q$. Since $u_q$ is a Hopf algebra, we have a
canonical homomorphism $H^\bu\to HH^\bu$, thus $D^\bu$ is an $H^\bu$ linear
category.
For an object $X\in D^\bu$ its {\em cohomological support}
$supp(X) \subset Spec(H^\bu)$ is the
set-theoretic support of the $H^\bu$ module $End^\bu(X)$.

\begin{Prop}\label{Dsupp} For an object $X\in D$ we have:
$$
X\in D_{\infty/2} \iff supp(X)\subseteq \n.$$
\end{Prop}

\proof
It is well known that $supp(X) \supset supp(Y)\cup supp(Z)$ provided that there
exists a distinguished triangle $Y\to X\to Z\to Y[1]$, thus the set of objects
with cohomological support contained in $\n$ forms
a full triangulated subcategory.
In view of Lemma \ref{fifi}, to check the implication $\Rightarrow$
it sufficies to check that $supp(X)\subset \n$ if
$X=CoInd_{u_q^+}^{u_q}(M)$ for some $M$.
 For such $X$ we have $Ext^\bu_{u_q}(X,X)
=Ext^\bu _{u_q^+}(X,M)$.
Moreover, it is not hard to check that this isomorphism is compatible
with the $H^\bu$ action, where the action on the right hand side
is obtained as the composition $H^\bu\to H^\bu(u_q^+)\to
Ext^\bu_{u_q^+}(M,M)$ and the canonical right action of
$Ext^\bu_{u_q^+}(M,M)$. Thus in this case  $Ext^\bu_{u_q}(X,X)$ is
set-theoretically supported on $\n$.

Assume now that $X\in D$ is such that $supp(X)\subseteq \n$.
To check that $X\in D_{\infty/2}$ it suffices to show that $Ext^\bu_{u_q^-}
(M,X)$ is finite dimensional for any $M\in D^b(u_q^+-mod)$.
It is a standard fact
%
%YULYA PEVTSOVA SAYS IT'S CLAIMED IN [GK] WITHOUT PROOF, BUT SHOULD FOLLOW
%FROM THERE ARGUMENT... IT MAY HAVE BEEN USED IMPLICITLY IN [ABG].
%
that $Ext^\bu(M_1,M_2)$ is a finitely generated $H^\bu(u_q^-)$ module
for any $M_1,M_2\in D^b(u_q^--mod)$,
thus it suffices to see that the $H^\bu(u_q^-)$-module
$Ext^\bu_{u_q^-}(M,X)$ is supported at $\{0\}\subset \n^-=Spec(H^\bu(u_q^-))$.
This is clear, since viewed as a $H^\bu(X)$ it is supported on $\n$.
\epf

\subsection{A description of the derived $u_q$-modules category
via coherent sheaves}

Let $\Nt=T^*(\calB)=\{(\b,x)\ |\ \b\in \calB, \, x\in rad(\b) \}$,
where $\calB=G/B$ is the flag variety of $G$ identified with the
set of Borel subalgebras in $\g$, and $rad$ stands for the nil-radical.
Let $\pi:\Nt\to \N$ be the Springer map, $\pi:(\b,x)\mapsto x$.

The result of \cite{Anichka} (based on \cite{ABG})
 yields a triangulated functor $\Psi:D^b (Coh^
{\Gm}(\Nt)) \to D^b(u_q-Mod)$, where $\Gm$ acts on $\Nt$ by
$t:(\b,x)\mapsto (\b,t^2x)$, and $u_q-Mod$ stands for the category
of finite dimensional modules. Notice that in contrast with
the definition of $u_q-mod$, the modules
in $u_q-Mod$  {\em do  not} carry a grading.\footnote{In fact,
$D^b (Coh^ {\Gm}(\Nt))$ can be identified with the derived
category of a block in the category of graded modules over $u_q$
compatible with a certain grading on $u_q$, defined
in \cite{AJS}. However, unlike the natural grading
by weights and its modifications,
this
grading is neither explicit, nor elementary; it is similar to
a grading on the category $O$ of $\g$ modules with highest weight arising
from Hodge weights on the Hom space between Hodge $D$-modules, or from
Frobenius weights.}

 The functor satisfies the following properties:
\begin{equation}\label{twist}
\Psi(\F(1))\cong \Psi(\F)[ 1] %{\rm\ \ \  check\  the\  sign?},
\end{equation}
where $\F(1)$ is the twist of $\F$ by the tautological character of $\Gm$;

\begin{equation}\label{equiv}
\Psi:\oplusl_{n\in \Zet} \Hom(\F, \G[n](n)) \iso \Hom (\Psi(\F), \Psi(\G));
\end{equation}

\begin{equation}\label{gen}
\langle Im(\Psi)\rangle = D^b(u_q-Mod_0),
\end{equation}
where $\langle Im(\Psi)\rangle$ denotes the full triangulated subcategory
generated by objects of the form $\Psi(\F)$, and $u_q-Mod_0$ is the
block (direct summand)
of the category $u_q-Mod$ which contains the trivial representation;

\begin{equation}\label{Otok}
\Psi(\O_{\Nt})=\k.
\end{equation}

The following slight generalization of this result is proved
by a straightforward modification of the argument of \cite{Anichka}.

\begin{Prop}
Let $C$ be a subtorus in the maximal torus $T$, and let $u_q-mod^C$
be the category of $u_q$-modules carrying a compatible grading by weights
of $C$. There exists a functor $\Psi^C:D^b(Coh^{C\times \Gm}(\Nt))
\to D^b(u_q-mod^C)$ satisfying properties \eqref{twist}-- \eqref{Otok} above.
\end{Prop}

\subsection{Semi-inifinite cohomology as cohomology with support}

From now on we fix $C$ to be a copy of the multiplicative group corresponding
to the coweight used to define the grading on $u_q$ (see the beginning
of this section), thus we have $u_q-mod^C=u_q-mod$.

\begin{Thm}\label{gamman}
For $\F\in D^b(Coh^{C\times \Gm})$ we have a canonical isomorphism
$$Ext^{\infty/2+i}_{u_q}(\k,\Psi(\F))\cong
R\Gamma_\n(\pi_*(\F)).$$
\end{Thm}

\proof We have
$$R\Gamma_\n^\bu (\pi_*(\F))\cong \varinjlim Ext^\bu (\O_\Nt/\pi^*(I),\F),$$
where $I$ runs over $C\times G_m$ invariant ideals in $\O_\N$
with support on $\n$. 
We have a canonical arrow $\Psi(\O_\Nt)\to \Psi(\O_\Nt/\pi^*(I))$, and
 in view of Proposition \ref{Dsupp} we have
$\Psi(\O_\Nt/\pi^*(I)) \subset D_{\infty/2}$. Thus by Theorem \ref{teorem}
we have a natural map 
%from the right hand side to the left hand side
%of the required isomorphism.
$$
R\Gamma_\n(\pi_*(\F)) \To Ext^{\infty/2+i}_{u_q}(\k,\Psi(\F)).
$$
In view of Proposition \ref{proind},
to check that this map is an isomorphism it suffices to show that
the pro-object $\widehat{\Psi(\O)}$ in $\Dgood$ defined by
$\widehat{\Psi(\O)}=\varprojlim \Psi(\O_\Nt/\pi^*(I)) $ represents
the same functor on $\Dgood$ as the object $\k=\Psi(\O)\in D$.

Let $X\in \Dgood$, and let $f_1, \dots, f_n$ be a regular sequence
in $\O(\N)$ whose common set of zeroes equals $\n$. We can and will assume
that $f_i$ is an eigen-function for the action of $C\times \Gm$.
There exists $N$ such that $f_i^N$ maps to $0\in End^\bu(X)$.
Then any morphism $\k\to X$ factors through $\k_N=\Psi(\O/(f_i^N))$.
This shows that the map $\varinjlim Hom (\k_N, X)\to Hom(\k,X)$ is 
surjective. Similarly, for large $N$ the map $Hom(\k_N,X)\to Hom (\k_{2N},X)$
kills the kernel of the map $Hom(\k_N,X)\to Hom (\k,X)$. Thus the map
  $\varinjlim Hom (\k_N, X)\to Hom(\k,X)$ is injective.
 \epf

\begin{Cor}
Let $T$ be a tilting module over Lusztig's ``big'' quantum group $U_q$.
The semi-infinite cohomology $Ext^{\infty/2+i}_{u_q}(\k,T)$ either
vanishes or is canonically isomorphic to $R\Gamma_\n(\F)$, where $\F\in
D^b(Coh^G(\N))$ is a certain (explicit) irreducible object 
in the heart of the perverse $t$-structure corresponding
to the middle perversity \cite{izvrat}.
\end{Cor}

\proof By the result of \cite{Humph} we have $T=\Psi(\Ftil)$
for some $\Ftil\in D^b(Coh^{G\times \Gm}
(\Nt))$, such that $\pi_*(\Ftil)$ either
vanishes or is an (explicit) irreducible perverse equivariant 
coherent sheaf as above.
The statement now follows from Theorem \ref{gamman}.\epf

\begin{Ex}
If $T=\k$ is the trivial module, then it is clear from the construction
of \cite{Humph} that we can set $\Ftil=\O_{\Nt}$. Thus
$\F\cong \O_\N$, so Corollary yields the main result
of \cite{Ar1}.
\end{Ex}

\end{section}

\end{document}